\documentclass[10pt]{article}

\usepackage{pifont}
\usepackage{amsmath}
\usepackage{cases}
\usepackage{amsfonts}
\usepackage{latexsym}
\usepackage{amssymb}
\usepackage{stmaryrd}
\usepackage{overpic}
\usepackage{color}
\usepackage{abstract}
\usepackage{indentfirst}
\usepackage{booktabs}
\usepackage{url}
\usepackage{verbatim}
\usepackage{diagbox}
\usepackage{multirow} 

\usepackage{graphicx}  
\usepackage{float}  
\usepackage{subfigure}

\newcommand{\bx}{\boldsymbol{x}}

\newcommand{\bH}{\boldsymbol{H}}
\newcommand{\bQ}{\boldsymbol{Q}}
\newcommand{\bp}{\boldsymbol{p}}
\newcommand{\bq}{\boldsymbol{q}}

\newcommand{\bn}{\boldsymbol{n}}

\newcommand{\ds}{\mathrm{d}s}

\DeclareFontEncoding{FMS}{}{}
\DeclareFontSubstitution{FMS}{futm}{m}{n}
\DeclareFontEncoding{FMX}{}{}
\DeclareFontSubstitution{FMX}{futm}{m}{n}
\DeclareSymbolFont{fouriersymbols}{FMS}{futm}{m}{n}
\DeclareSymbolFont{fourierlargesymbols}{FMX}{futm}{m}{n}
\DeclareMathDelimiter{\VERT}{\mathord}{fouriersymbols}{152}{fourierlargesymbols}{147}

\newtheorem{theorem}{Theorem}[section]

\newtheorem{example}[theorem]{Example}

\newtheorem{remark}[theorem]{Remark}

\begin{document}

\begin{center}
\large\bf Deep Petrov-Galerkin Method for Solving Partial Differential Equations
\end{center}

\vspace{2mm}

\begin{center}
 {\large\sc Yong Shang}\footnote{School of Mathematics and Statistics, Xi'an Jiaotong University, Xi'an,
Shaanxi 710049, P.R. China. E-mail: {\tt fsy2503.xjtu@xjtu.edu.cn}},\quad
{\large\sc Fei Wang}\footnote{School of Mathematics and Statistics, Xi'an Jiaotong University,
Xi'an, Shaanxi 710049, China. The work of this author was partially supported by
the National Natural Science Foundation of China (Grant No.\ 12171383). Email: {\tt feiwang.xjtu@xjtu.edu.cn}}\quad {\rm and}\quad
 {\large\sc Jingbo Sun}\footnote{School of Mathematics and Statistics, Xi'an Jiaotong University, Xi'an,
Shaanxi 710049, P.R. China. E-mail: {\tt jingbosun.xjtu@xjtu.edu.cn}}
\end{center}

\vspace{2mm}

\begin{quote} 
\noindent{}{\bf Abstract}:
Deep neural networks are powerful tools for approximating functions, and they are applied to successfully solve various problems in many fields. In this paper, we propose a neural network-based numerical method to solve partial differential equations. In this new framework, the method is designed on weak formulations, and the unknown functions are approximated by deep neural networks and test functions can be chosen by different approaches, for instance, basis functions of finite element methods, neural networks, and so on. Because the spaces of trial function and test function are different, we name this new approach by Deep Petrov-Galerkin Method (DPGM). The resulted linear system is not necessarily to be symmetric and square, so the discretized problem is solved by a least-square method. Take the Poisson problem as an example, mixed DPGMs based on several mixed formulations are proposed and studied as well. In addition, we apply the DPGM to solve two classical time-dependent problems based on the space-time approach, that is, the unknown function is approximated by a neural network, in which temporal variable and spatial variables are treated equally, and the initial conditions are regarded as boundary conditions for the space-time domain. Finally, several numerical examples are presented to show the performance of the DPGMs, and we observe that this new method outperforms traditional numerical methods in several aspects: compared to the finite element method and finite difference method, DPGM is much more accurate with respect to degrees of freedom; this method is mesh-free, and can be implemented easily; mixed DPGM has good flexibility to handle different boundary conditions; DPGM can solve the time-dependent problems by the space-time approach naturally and efficiently. The proposed deep Petrov-Galerkin method shows strong potential in the field of numerical methods for partial differential equations.


\end{quote}

\vspace{2mm}

\section{Introduction}

Deep neural networks (DNNs) have achieved great success in artificial intelligence, scientific computing and other machine learning tasks (\cite{hinton2012ml7,lecun2015ml2,schmidhuber2015ml3,goodfellow2016ml1,razzak2018ml4,voulodimos2018ml5,young2018ml6,zhang2018ml8}). DNNs produce a large class of nonlinear functions through compositional construction. Due to their powerful universal approximation ability, in recent years, DNNs have been applied for solving partial differential equations (PDEs), and several DNN-based methods (\cite{ee2017spde1,berg2018spde4,han2018spde2,he2018spde6,ee2018drm,raissi2019pinn,chen2019spde5,fan2019spde7,khoo2019spde3,long2019spde10,li2020spde8,liu2020spde9}) were proposed to overcome the difficulty so-called the ``curse of dimensionality" of the traditional PDE solvers such as finite element method (FEM), which requires a discretization of the interested domain, while the number of the mesh points will increase exponentially fast with respect to the problem dimension and make it quickly become computationally intractable. In such a situation, the generation of meshes is very time-consuming as well.

Among these DNN-based methods, one approach is to establish the algorithm with the strong form and minimize the residual of the original PDEs, such as physical informed neural networks (PINNs), deep Galerkin method (DGM), mixed residual method (MIM), local extreme learning machines (locELM) and so on. Specifically, PINN was proposed in \cite{raissi2019pinn} to use neural networks to approximate the solution of PDEs, and the neural networks are trained by minimizing the mean squared errors of the residual of the observed data, the differential equation, boundary conditions, and/or initial conditions. After that, some other model variants, like fractional PINNs and nonlocal PINNs are put forward in \cite{pang2019fpinn,pang2020npinn,zhang2020tpinn}. DGM designs the loss function as PDE residual in the least-square sense to measure how well the approximate solution satisfies the differential operator and boundary conditions (\cite{sirignano2018dgm}). MIM rewrites the PDE into a lower-order system and then uses the system residual as the loss function (\cite{lyu2022mim}). Combining the ideas of domain decomposition and extreme learning machines (ELM), locELM achieves high-precision results by solving the parameters of the neural network with least-square computations rather than training by an optimizer (\cite{dong2021elm}).

Another approach focuses on the weak formulations of the PDEs. Based on the variational principle, the deep Ritz method (DRM) uses the variational formulation of symmetric PDEs and integral of the boundary residual as the loss function, and the DNNs are trained by minimizing the numerical quadrature of the loss function (\cite{ee2018drm}). The deep Nitsche method (\cite{liao2019dnm}) adopts Nitsche's formulation as the loss function to deal with the essential boundary conditions rather than a penalty in DRM. Based on the Euler-Lagrange equation of the energy functional, a penalty-free neural network method (\cite{sheng2021pfnn}) is proposed to solve the second-order boundary-value problems on complex geometries. Weak adversarial networks (\cite{zang2020wan}) reformulate the PDEs as a saddle-point problem in the weak formulation and approximate trial and test functions by two neural networks and then train them alternately as an operator norm minimization problem.

These DNN-based methods have made achievements in some aspects, such as being able to solve problems in high dimensions or on irregular domains, including the observed data into the loss function so that it can handle inverse problems easily. However, some difficulties arise and need to be explored further. One of the biggest problems is the accuracy limitation of the numerical solutions obtained by these methods due to the weakness of the optimizer for the training process, even DNNs have very nice approximation properties (\cite{jagtap2020limit}). Meanwhile, the training process could cost a huge amount of time. For example, a DNN-based method may take hours to train the neural network to reach a certain accuracy, and it just takes a few seconds for FEM. They are not sufficiently efficient in solving general PDEs, and more researches need to be explored further.

On the other hand, providing rigorous error analysis like FEM is a very challenging task for these DNN-based methods. For PINNs, the convergence analysis is provided in \cite{mishra2020errorpinn3,shin2020errorpinn1, shin2020errorpinn2}, and a PINN with ReLu$^3$ network is analyzed and the convergence rate was given in $C^2$ norm (\cite{jiao2021errorpinn4}). The error analysis of DRM was established in \cite{luo2006barron1,xu2020barron4,hong2021barron3,lu2021barron2} via assuming that the exact solution is contained in the spectral Barron space which has the property of being approximated by a two-layer neural network. The convergence rate of DRM with smooth activation functions like logistic or hyperbolic tangent was derived in $H^1$ norm for elliptic equations (\cite{duan2021errordrm1,jiao2021errordrm2}). Such analysis can help us comprehend the convergence property with respect to the depth, width, and parameters in the neural networks.

Least-square methods have been studied for solving partial differential equations as well (\cite{aziz1985lsq1,carey1988lsq5,bochev2001lsq2,bochev2001lsq4,sterck2004lsq6,bochev2009lsfem,bochev2016lsq3}). By employing the least-square functional as its loss function, a least-square neural network method (LSNN) is proposed to solve the advection-reaction problem with discontinuous solution (\cite{cai2021lsnn}), and deep least-square method based on the least-square functional of the first-order system is studied in \cite{cai2020dlsq}. It is worth mentioning that the locELM can be regarded as a least-square method by obtaining the parameters of the neural network through solving a least-square problem, and it appears to be more accurate and computationally efficient (\cite{dong2021elm}). 

In this paper, based on variational principle, we proposed a Deep Petro-Galerkin Method (DPGM) in which the numerical solutions are approximated by DNNs while test functions can be chosen by different approaches. The parameters of the neural network are given randomly and fixed except for the last layer, which not only reduces the parameters that need to be trained but also facilitates the assembling process of the linear system, then we solve it by a least-square method. This approach indeed improves the accuracy of the numerical solution and reduces the computational cost as well. Meanwhile, for solving time-dependent problems, temporal and spatial variables are treated jointly and equally, and we adopt the space-time approach in the DPGM framework so that the initial conditions are treated as boundary conditions. Therefore, DPGM can solve time-dependent problems without any iteration steps, and the resulted neural networks can supply the value of the numerical solution at any given space-time point without interpolation as the traditional numerical methods do.

To deal with the boundary conditions, one approach is the penalty method adopted by PINNs, DRM, DGM, and some other methods, and another way is to construct a special neural network satisfying the boundary condition intentionally (\cite{lagaris1988bc3,lagaris2000bc4,lyu2008mimb,berg2018spde4,khoo2019spde3}). However, the first approach is suffering difficulty of choice of penalty parameters, which is crucial for balancing the residual of the interior and boundary items; and the second one is not practical for the domain with complicated boundary geometry. Instead, in DPGM, we adopt the ideas of the least-square method by assembling the boundary conditions and variational formulation to the linear systems simultaneously and solving it by a least-square method to determine the unknown parameters. In this way, no penalty parameter is involved and the boundary conditions can be satisfied easily. Furthermore, DPGM maintains the advantages of the DNN-based methods, for example, the computation for the numerical quadrature can be adopted by quasi-Monte Carlo method (\cite{niederreiter1992qmc1,caflisch1998qmc2,dick2013qmc3,chen2019spde5}), which produces a mesh-free method and is essential for high-dimensional PDEs, and thus DPGM can be used to solve the PDEs on complicated geometries. The proposed DPGMs have the following advantages: neural networks are used to approximate the numerical solution to guarantee excellent approximation property; the ideas of extreme learning machine and least-square method are adopted so that the resulted discrete problem can be solved efficiently; both the essential and natural boundary conditions can be treated easily; DPGM can be applied with mixed formulations without worrying the stability issues through solving the discrete problem by the least-square method; time-dependent problems can be solved by DPGM in space-time approach efficiently.

The rest of this paper is organized as follows. In Section 2, we introduce the basic ideas of the DPGM for solving an elliptic partial differential equation. In Section 3, taking the Poison problem as an example, we establish mixed-DPGM with four mixed formulations. Then we consider solving time-dependent problems by DPGM in the space-time approach in Section 4. In Section 5, numerical examples are provided to show the efficiency and accuracy of the proposed DPGM. Finally, conclusions and discussions are drawn in the last section.

\section{Deep Petrov-Galerkin Method}\label{prem}

In this section, we introduce the basic ideas of the deep Petrov-Galerkin method and show how to apply it to solve an elliptic partial differential equation.

Consider a partial differential equation 

\begin{align}
	&\mathcal{A} u = f \qquad {\rm in}\; \Omega,\label{pde1} \\
	&\mathcal{B} u = g \qquad{\rm on}\; \Gamma, \label{pde2}
\end{align}
where $\Omega$ is a bounded domain in  $\mathbb{R}^d$, and its boundary $\Gamma$ is split as $\Gamma = \Gamma_D \cup \Gamma_N $ with ${\Gamma}_D \cap {\Gamma}_N = \emptyset$.
For example, when the differential operator $\mathcal{A}$ is given by
\begin{equation}\label{diffreac}
	\mathcal{A} := -\nabla\cdot(\alpha(\bx)\nabla ) + \delta(\bx),
\end{equation}
where $\alpha_0\geq\alpha(\bx)\geq \alpha_0>0$, $\delta_1 \geq \delta(\bx) \geq \delta_0>0$, \eqref{pde1} is a 
diffusion-reaction equation. Here, $\alpha_0$, $\alpha_1$, $\delta_0$ and $\delta_1$ are constants.
The boundary condition \eqref{pde2} can be Dirichlet, Neumann and Robin types. In this paper, we consider the following mixed boundary conditions:
\begin{align}
	u = g_D& \qquad {\rm on}\; \Gamma_D,\label{bc1} \\
	\alpha(\bx)\nabla u\cdot \bn = g_N &\qquad {\rm on}\; \Gamma_N. \label{bc2}
\end{align}

\subsection{Petrov-Galerkin method}

Under proper conditions, the PDE \eqref{pde1} with boundary conditions \eqref{bc1}--\eqref{bc2} has the week formulation:
Find $u\in U_{D,g_D} = \{u\in U;\; u|_{\Gamma_D} = g_D\}$ such that
\begin{align}\label{prb1}
	a(u,v) = l(v)\qquad \forall v\in V_{D,0} = \{v\in V;\; v|_{\Gamma_D} = 0\},
\end{align}
where $U$ and $V$ are two Hilbert spaces with inner products $(\cdot,\cdot)_U$ and $(\cdot,\cdot)_V$, as well as corresponding norms  
$\|\cdot\|_U$ and $\|\cdot\|_V$, respectively. Here, $a(\cdot,\cdot): U\times V\rightarrow \mathbb{R}$ is a continuous bilinear form with 
\begin{equation}
	a(u,v) \leq M \Vert u\Vert_{U}\Vert v\Vert_{V}  \qquad \forall u\in U,\, v\in V\label{assum1}
\end{equation}
for a constant $M>0$. By Babu\v{s}ka Theorem (\cite{babuska1972,babuska1973}), the problem \eqref{prb1} has a unique solution if and only if the following conditions hold
\begin{align}
	\inf\limits_{u\in U}\sup\limits_{v\in V} \frac{a(u,v)}{\Vert u\Vert_{U}\Vert v\Vert_{V}} =
	\inf\limits_{v\in V}\sup\limits_{u\in U} \frac{a(u,v)}{\Vert u\Vert_{U}\Vert v\Vert_{V}} =\beta >0, \label{assum2}
\end{align}
which is known as the Babu\v{s}ka--Brezzi or inf-sup condition (\cite{boffi2013mixed}). If we choose $\mathcal{A}$ as \eqref{diffreac}, the bilinear form and linear form in \eqref{prb1} are given by
\begin{align*}
a(u,v)&=\int_\Omega \left(\alpha(\bx)\nabla u\cdot \nabla v + \delta(\bx) u\, v\right)\,{\rm d}\bx,\\
l(v) &= \int_\Omega f\, v\,{\rm d}\bx + \int_{\Gamma_N} g_N\, v\,\ds,
\end{align*}
and the bilinear form $a(u,v)$ satisfies the conditions \eqref{assum1} and \eqref{assum2}.

To numerically solve the problem \eqref{prb1}, one can apply the Petrov-Galerkin (PG) method:
Find $u_h\in U_h\subset U_{D,g_D}$ such that
\begin{align}
	a(u_h,v_h) = l(v_h)\qquad \forall v_h\in V_h\subset V_{D,0}.\label{prb2}
\end{align}
Here, $U_h$ and $V_h$ are finite-dimensional function spaces, and usually consist of piecewise polynomial functions. By \eqref{assum2}, we know that the problem \eqref{prb2} is well-posed if and only if the following conditions hold
\begin{align}
	\inf\limits_{u_{h}\in U_{h}}\sup\limits_{v_h\in V_h} \frac{a(u_{h},v_h)}{\Vert u_{h}\Vert_{U}\Vert v_h\Vert_{V}} =  
	\inf\limits_{v_h\in V_h}\sup\limits_{u_{h}\in U_{h}} \frac{a(u_{h},v_h)}{\Vert u_{h}\Vert_{U}\Vert v_h\Vert_{V}} =  \beta_h >0. \label{assum4}
\end{align}
Furthermore, when $U_h$ and $V_h$ are finite-dimensional, the above two conditions are reduced to one (\cite{xu2003brezzi}). Then the fundamental result for the PG method is given as follows.
\begin{theorem}[\cite{babuska1972}]
	If the conditions \eqref{assum2} and \eqref{assum4} hold, and $u$ and $u_{h}$ are the solutions of the problems \eqref{prb1} and \eqref{prb2}, respectively, then 
	\begin{equation}
		\Vert u-u_{h}\Vert_{U} \leq (1+\frac{M}{\beta_h}) \inf\limits_{w_{h}\in U_{h}} \Vert u-w_{h}\Vert_{U}.
	\end{equation}
\end{theorem}

\begin{remark}
Note that the finite element pair of $U_h$ and $V_h$ needs to be elaborately designed for the PG method, and the inf-sup condition \eqref{assum4} has to be verified carefully, which are very challenging and technical tasks for general users. Can we avoid these difficulties? Let us try to do something different in this paper.
\end{remark}

\subsection{Neural networks}
\label{NNs}

There are a variety of neural network structures, let us introduce the neural networks used in this paper as follows. 
Let $D,n_i\in \mathbb{N}^{+}$, and $n_i$ is the number of neurons in the $i$-th layer, $\rho$ is the activation function. 
A fully connected neural network with the depth of $D$ is a function $\mathbf{\Phi}:\mathbb{R}^{n_0}\rightarrow \mathbb{R}^{n_D}$ defined by
\begin{align}
	&\mathbf{\Phi}_0(\bx) = \bx,\nonumber\\
	&\mathbf{\Phi}_{l}(\bx) = \rho(\mathbf{W}_l\mathbf{\Phi}_{l-1}+\mathbf{b}_l) \qquad \text{for}\ l=1,\cdots,D-1,\label{neural}\\
	&\mathbf{\Phi}:= \mathbf{\Phi}_D(\bx) =\mathbf{W}_D\mathbf{\Phi}_{D-1} + \mathbf{b}_D,\nonumber
\end{align}
where $\mathbf{W}_l = \left(w_{ij}^{(l)} \right)\in \mathbb{R}^{n_l \times n_{l-1}}$ and $\mathbf{b}_l = \left(b_i^{(l)}\right)\in \mathbb{R}^{n_l}$, $\phi_l = \left\{\mathbf{W}_l,\mathbf{b}_l\right\}$ are called the weight parameters in $l$-th layer. In this paper, we set $n_0 = d$ for stationary problems and $n_0 = d + 1$ for time-dependent problems.

In addition, we also adopt the ResNet (\cite{he2016resnet}), in which the input layer is a fully-connected layer with the number of neurons $n_0=d$, each layer of this network is constructed by stacking several blocks, each block consists of two linear transformations, two activation functions, and a residual connection. Thus a ResNet $\mathbf{\Psi}$ with $D$ layers takes the form as:
\begin{align*}
\qquad\qquad\qquad\qquad\qquad\;\;\;\;	&\mathbf{\Psi}_0(\bx) = \rho(\mathbf{W}_0\mathbf{\bx}+\mathbf{b}_0) ,\\
	&\mathbf{\Psi}_{l}(\bx) =\rho\left (\mathbf{W}_{l,2}\,\rho(\mathbf{W}_{l,1}\mathbf{\Psi}_{l-1}+\mathbf{b}_{l,1})+\mathbf{b}_{l,2}\right) +\mathbf{\Psi}_{l-1} \qquad \text{for}\ l=1,\cdots,D-1,\\
	&\mathbf{\Psi}:= \mathbf{\Psi}_D(\bx) =\mathbf{W}_D\mathbf{\Psi}_{D-1} + \mathbf{b}_D,
\end{align*}
where $\{\mathbf{W}_{l,1},\mathbf{b}_{l,1}\}$ and $\{\mathbf{W}_{l,2},\mathbf{b}_{l,2}\}$ are the weight parameters of the first and second linear transformation in the $l$-th layer, respectively.

Denote the number of nonzero entries of a matrix $\mathbf{W}_k$ by 
$$\Vert \mathbf{W}_k \Vert_{l^0}:= \left|\{(i,j): w_{ij}^{(k)}\neq 0\}\right|,$$
then $$N_D:= \sum\limits_{j=1}^{D}(\Vert \mathbf{W}_j \Vert_{l^0} + \Vert \mathbf{b}_j \Vert_{l^0})$$
is the total number of nonzeros weights of the neural network $\mathbf{\Phi}$. For two positive constant integers $M_D$ and $B_{D}$, we give a class of neural network functions as follows
\begin{align*}
	\mathcal{N}_{\rho}(D,M_D,B_{D}) :=\{ \mathbf{\Phi}\; \text{defined by \eqref{neural} with depth D}, N_D\leq M_D\;  \text{and}\;
	|w_{ij}^{(l)}|\leq B_D, |b_i^{(l)}|\leq B_D \}.
\end{align*}

The universal approximation theorem (\cite{cybenko1989appro1,hornik1991appro2}) clarifies that every continuous function on a compact domain can be uniformly approximated by shallow neural networks with continuous, non-polynomial activation functions. The relationship between ReLU-DNN and linear finite element function was studied in \cite{he2018spde6}.
More results have been established in \cite{barron1994acti1, mhaskar1996acti3,shaham2018acti4,bolcskei2019acti2,tang2019acti5} for activation functions with a certain regularity, and these approximation errors were given in the sense of $L^p$ norm. Then the error bounds for H{\"o}lder functions and functions in the Sobolev space $W^{n,\infty}$ were given under $L^{\infty}$ norm (\cite{yarotsky2017wp,ohn2019holder}). Finally, the error bounds for ReLU neural networks and ELU-neural networks are derived in Sobolev norms (\cite{guhring2020wkp1,guhring2021wkp2}), and it holds for many practically used activation functions such as the logistic function, tanh, arctan, and others, which offers us a powerful tool to study the error analysis for those DNN-based methods.

\begin{theorem}[Proposition 4.8 in \cite{guhring2021wkp2}, Theorem 4.1 in \cite{jiao2021errordrm2}] \label{3.0}
	Given $p\geq 1$, $s,k,d\in \mathbb{N}^{+}$, $s\geq k+1$. Let $\rho$ be the logistic function $\frac{1}{1+e^{-x}}$ or tanh function $\frac{e^x-e^{-x}}{e^x+e^{-x}}$. For any $\epsilon >0$ and $f\in \mathcal{F}_{s,p,d}$, there exists a neural
	network $f_{\rho} \in \mathcal{N}_{\rho}(D,M_D,B_{D})$  with depth $D\leq C\,{\rm log}(d+s)$, $M_D\leq C\cdot\epsilon^{-d/s-k-\mu k}$, and $B_D\leq C\cdot {\epsilon}^{-\theta}$ such that
	\begin{equation*}
		\Vert f-f_{\rho} \Vert_{W^{k,p}([0,1]^d)} \leq \epsilon.
	\end{equation*}
	where $C,\theta$ are constants depending on d,s,p,k; $\mu$ is an arbitrarily small positive number and $\mathcal{F}_{s,p,d}:=\{f\in W^{s,p}([0,1]^d):\Vert f\Vert_{W^{s,p}([0,1]^d)} \leq 1\}$.
\end{theorem}

\begin{remark}
	The  bounds in the above theorem can be found in the proof of Proposition 4.8 in \cite{guhring2021wkp2}, and the bound on the depth was given in Theorem 4.1 in \cite{jiao2021errordrm2} explicitly. With these results, we see that the solutions of PDEs can be approximated very well by neural networks with a sufficient number of layers and nonzero weights.
\end{remark}

\subsection{Deep Petrov-Galerkin method}
\label{NNs}

Instead of approximating the exact solution $u(\bx)$ by a piecewise polynomial function $u_h$, in this paper,
we try to search a neural network $u_\rho(\bx;\theta)$ to approximate $u(\bx)$ with proper parameters $\theta$. Let us propose a numerical scheme: Find $u_\rho \in U_{\rho} \subset \mathcal{N}_{\rho}(D,M_D,B_{D}) \subset U$ such that 
\begin{align}
	a(u_\rho,v_h) &= l(v_h)\qquad\; \forall v_h\in V_h\subset V_{D,0},\label{prb3}  \\ 
	u_\rho(\bx_k) &= g_D(\bx_k)  \quad\; {\rm for\ some\ points}\ \bx_k\in \Gamma_D,\; k=1,2,\cdots, N_b,\label{boundary_D}  
\end{align}
which is named by Deep Petrov-Galerkin Method (DPGM). Here, the trial function space $U_{\rho}$ refers to a class of neural network functions with the activation function $\rho$, and the structure is introduced in the previous subsection. 
In addition, the test function space $V_{h}$ can be chosen as any proper finite-dimensional function space, which should be a good approximation of $V_{D,0}$. In this paper, we choose $V_h$ as finite element spaces, and mainly use the bilinear finite element space for 2D problems and trilinear finite element space for 3D problems. Without loss of generality, we assume that $\{v_j\}_{j=1}^G$ is a group of functions belonging to $V_{D,0}$ and $\bigcup_{j=1}^G {\rm supp}\{v_j\} =\Omega$, then we choose
$V_h:={\rm span}\{v_j\}_{j=1}^G$.

\begin{figure}[!htbp] 		
	\centering
	\includegraphics[scale=0.55]{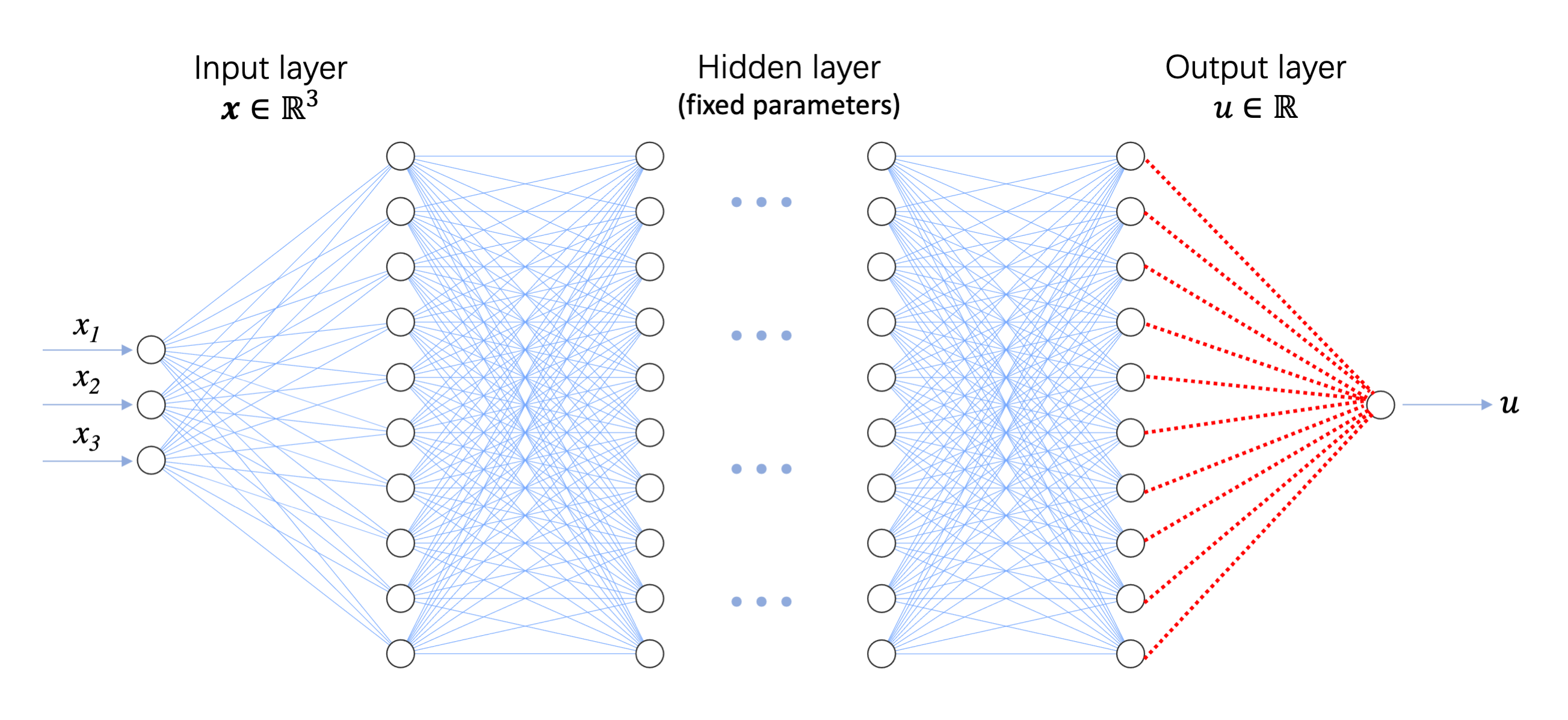}
	\caption{Network structure of $u_\rho:\mathbb{R}^{3}\rightarrow \mathbb{R}$, the solid blue line represents the parameters of the neural network that are randomly initialized and fixed afterwards, the dotted red line refers to the parameters that need to be solved. }
	\label{network}
\end{figure}

Notice that $\mathcal{N}_{\rho}(D,M_D,B_{D})$ is not a function space in general, so we did not choose it as $U_\rho$, instead, we choose a function space, a subset of $\mathcal{N}_{\rho}(D,M_D,B_{D})$ as $U_\rho$, which is constructed as follows and showed in Figure \ref{network}. For any $u_\rho\in U_\rho\subset \mathcal{N}_{\rho}(D,M_D,B_{D})$, we choose the initial weights $\phi_l$ with values drawn from the uniform distribution $\mathcal{U}(-r, r)$ for each layer and fixed them except for the last layer, where $r\in \mathbb{R}$. Let $\Phi_{D-1}^j$ denote the output of the $j$-th neuron in the $D-1$ layer ($j=1,2,\cdots,n_{D-1}$),  $\mathbf{W}_D = \left((w_{ij}^{(D)}) \right)\in \mathbb{R}^{n_D\times n_{D-1}}$ with $n_D = 1$, and $\mathbf{b}_D=\mathbf{0}$. Then the output of the neural network 
\begin{equation*}
	\mathbf{\Phi}:= \mathbf{\Phi}_D(\bx) =\mathbf{W}_D\mathbf{\Phi}_{D-1} + \mathbf{b}_D =  \sum\limits_{j=1}^{ n_{D-1}}
	w_{1j}^{D} \Phi_{D-1}^j(\bx).
\end{equation*}
For simplicity, we denote $\Phi_j^u :=\Phi_{D-1}^j$ and $u_\rho^j :=w_{1j}^{D}$ for $j = 1,\cdots,n_{D-1}$. Then 
\begin{equation}
	u_\rho(\bx) = \sum\limits_{j=1}^{ n_{D-1}} u_\rho^j \Phi_j^u(\bx). \label{nnappu}
\end{equation}
Notice that the values of the initial weights $\phi_l$ ($l = 1,\cdots,n_{D-1}$) are given randomly and fixed, so these weights do not need to be trained or solved, then $\{\Phi_j^u(\bx)\}_{j=1}^{n_{D-1}}$ can be regarded as a group of basis functions. Hence,
$$U_\rho = {\rm span}\{\Phi_1^u(\bx),\cdots,\Phi_{n_{D-1}}^u(\bx)\}.$$
We see that the degrees of freedom for the DPGM is the number of neurons in $D-1$ layer, that is, $n_{D-1}$.

Then, the problem \eqref{prb3} becomes: Find $u_\rho^j $, $j = 1,\cdots,n_{D-1}$, such that
\begin{equation}
	\sum\limits_{j=1}^{ n_{D-1}}u_\rho^j a(\Phi_j^u(\bx),v_i) = l(v_i)\qquad \forall v_i\in V_h\subset V_{D,0}.  \label{linear1}
\end{equation}
The above equation gives 
$$\mathbf{A} \boldsymbol{U}=\mathbf{L},$$ 
where $\mathbf{A}_{i,j} = a(\Phi_j^u(\bx),v_i)$, $\boldsymbol{U}= (u_\rho^1,\cdots,u_\rho^{n_{D-1}})^{T}$, and $\mathbf{L}=(l_1,\cdots,l_{N_h})^{T}$ with $l_i = l(v_i)$ ($i=1,\cdots,N_h$). Note that the Neumann boundary condition \eqref{bc2} is included in the weak formulation. However, we need to enforce the Dirichlet boundary condition \eqref{boundary_D}. To do so, we take some random samples $\{\bx_k\}_{k=1}^{N_b}$ according to the uniform distribution $\mathcal{U}(\Gamma_D)$ , and set $u_\rho(\bx_k)  =g_D(\bx_k)$ for all $ \bx_k \in \Gamma_D$, i.e.,
\begin{equation}
	\sum\limits_{j=1}^{ n_{D-1}}u_\rho^j \Phi_j^u(\bx_k)= g_D(\bx_k)  \qquad \text{for} \  k = 1,\cdots,N_b. \label{linear2}
\end{equation}
Similarly, the above equation implies $\mathbf{B} \boldsymbol{U}=\mathbf{G}$, where $\mathbf{B}_{k,j} = \Phi_j^u(\bx_k)$ and $\mathbf{G}=(g_1,\cdots,g_{N_b})^{T}$ with  $g_k  =g_D(\bx_k)$.
Thus, we can obtain the solution $u_\rho$ by solving a least-square problem with the linear system 
\begin{equation}
\left[\begin{array}{c}\mathbf{A} \\\mathbf{B} \end{array}\right] \boldsymbol{U}
=\left[\begin{array}{c}\mathbf{L} \\\mathbf{G} \end{array}\right].\label{linear3}
\end{equation}

The algorithm for the DPGM  is summarized as follows.

\begin{tabular*}{16cm}{l}
	\hline
	{\bf  Algorithm 1} Deep Petrov-Galerkin Method \\ \hline
	\quad Step 1. Initialize network architecture $u_{\rho} :\Omega \rightarrow \mathbb{R}$ with depth $D$ and parameters ${\phi_l}$ for $ l=1,\cdots,D$. \\
	\quad Step 2. Fix the parameters ${\phi_l}$ for $ l=1,\cdots,D-1$ of $u_{\rho}$ and rewrite  $u_{\rho}$ as \eqref{nnappu}.\\
	\quad Step 3. Choose $v_i \in V_h$, $i = 1,\cdots,N_h$ and assemble the linear system  $\mathbf{A} \boldsymbol{U}=\mathbf{L}$ w.r.t.  \eqref{linear1}.\\
	\quad Step 4. Take random samples $\{\bx_k\}_{k=1}^{N_b}$  according to $\mathcal{U}(\Gamma_D)$ and obatin $\mathbf{B} \boldsymbol{U}=\mathbf{G}$  w.r.t. \eqref{linear2}.\\
	\quad Step 5. Obatin $\boldsymbol{U}$ by solving the least-square problem w.r.t.
	\eqref{linear3}.\\
	\quad Step 6. Update the network parameters $\phi_D$ of $u_{\rho}$.\\ \hline
\end{tabular*}

\begin{remark}
Most of the DNN-based numerical methods train neural networks by solving an optimization problem even the original problem may be a linear PDE, so the training process cost a large amount of computation time, and the precision of the results is less satisfied because it is very hard to find good solver for such optimization problems. If we want to avoid such a situation and also use powerful neural networks to approximate the solutions of PDEs, then the idea of ELM can be adopted so that we only need to solve a linear system. Note that the matrix $\mathbf{A}$ is not square and symmetric, so we cannot solve the discrete problem by usual linear system solvers. In addition, the essential boundary condition is not easy to be built into the structure of the neural networks, and the penalty method of enforcing the boundary condition may bring unexpected errors. By using the least-square method, we can get around these problems, and solve the linear system. The numerical experiments display that the hybrid of ELM and least-square approaches in Petrov-Galerkin formulation produces a highly efficient numerical method.
\end{remark}

\section{Mixed DPGM}

Similar to the mixed finite element method (\cite{raviart1977mixed}), we can construct a Mixed DPGM (M-DPGM) by approximating more than one unknown functions by one or more neural networks, simultaneously.

Take the Poisson equation with the mixed boundary conditions as an example
\begin{numcases}{}
-\nabla \cdot(\nabla u) = f \qquad\;\;\, \rm{in} \; \Omega \label{prb6} , \\
\qquad\qquad	u =g_D \qquad \rm{on} \; \Gamma_D \label{bc3},\\
\quad\;\;	\nabla u\cdot \bn\, = g_N\quad\;\;\; \rm{on} \; \Gamma_N. \label{bc4}
\end{numcases}
Set $\bp = \nabla u$, we have the following first-order system
\begin{numcases}{}
	\bp - \nabla u = 0 \qquad \rm{in} \, \Omega, \label{prb7}\\
	\; -\nabla \cdot \bp = f \quad\;\;\ \rm{in} \; \Omega, \label{prb8} \\
	\quad\;\;\;\;\;\; u =g_D \;\;\; \,\,\,\rm{on} \; \Gamma_D, \label{bc5}\\
	\quad\bp\cdot \bn\,= g_N \;\;\;\;\; \rm{on} \; \Gamma_N. \label{bc6}
\end{numcases}
On both sides of the equations \eqref{prb7}--\eqref{prb8}, multiplying proper test functions, we can obtain the following mixed formulations by integration by parts and the boundary conditions \eqref{bc5}--\eqref{bc6}. 
 
{\bf Mixed Formulation 1}:  
Find $(\bp,u) \in (L^2(\Omega))^d\times H^1_{D,g_D}(\Omega)$ such that 
\begin{align*}
	\int_{\Omega}\bp\cdot \bq\ {\rm d}\bx - \int_{\Omega} \nabla u\cdot \bq\ {\rm d}\bx &= 0 \qquad\qquad\qquad\qquad\quad\qquad \forall \bq\in (L^2(\Omega))^d, \\
	\int_{\Omega} \bp\cdot \nabla v \ {\rm d}\bx  &= \int_{\Omega} f\,v \ {\rm d}\bx +\int_{\Gamma_N}g_N\, v \ \ds \quad \forall v \in H^1_{D,0}(\Omega).
\end{align*}

{\bf Mixed Formulation 2}:  
Find $(\bp,u) \in \bH_{N,g_N}(\mathrm{div},\Omega)\times L^2(\Omega)$ such that 
\begin{align*}
	\int_{\Omega}\bp\cdot \bq \ {\rm d}\bx + \int_{\Omega}u \,\nabla \cdot \bq \ {\rm d}\bx&= \int_{\Gamma_D}g_D\,\bq\cdot \bn\ \ds \quad \forall \bq\in \bH_{N,0}(\mathrm{div},\Omega), \\
	-\int_{\Omega} \nabla \cdot \bp\, v\ {\rm d}\bx &= \int_{\Omega} f\,v \ {\rm d}\bx \qquad\quad\;\; \forall v \in L^2(\Omega).
\end{align*}

{\bf Mixed Formulation 3}:  Find $(\bp,u) \in \bH_{N,g_N}(\mathrm{div},\Omega)\times H^1_{D,g_D}(\Omega)$ such that 
\begin{align*}
	\int_{\Omega}\bp\cdot \bq\ {\rm d}\bx - \int_{\Omega} \nabla u\cdot \bq \ {\rm d}\bx&= 0 \qquad\qquad\;\; \forall \bq\in (L^2(\Omega))^d, \\ 
	-\int_{\Omega} \nabla \cdot \bp \,v \ {\rm d}\bx&= \int_{\Omega} f \,v \ {\rm d}\bx\quad \forall v \in L^2(\Omega).
\end{align*}

{\bf Mixed Formulation 4}:  Find $(\bp,u) \in (L^2(\Omega))^d\times L^2(\Omega)$ such that
\begin{align*}
\int_{\Omega}\bp\cdot \bq \ {\rm d}\bx+ \int_{\Omega}u \,\nabla \cdot \bq \ {\rm d}\bx&= \int_{\Gamma_D}g_D\, \bq\cdot \bn\ \ds \qquad\qquad\;\;\; \forall \bq\in \bH_{N,0}(\mathrm{div},\Omega), \\ 
\int_{\Omega} \bp\cdot \nabla v \ {\rm d}\bx  &= \int_{\Omega} f\,v \ {\rm d}\bx + \int_{\Gamma_N}g_N\,v \ \ds \quad \forall v \in H^1_{D,0}(\Omega).
\end{align*}
Here,
\begin{align*}
H^1_{D,0}(\Omega) &= \{v\in H^1(\Omega);\; v|_{\Gamma_D} = 0\},\\ 
H^1_{D,g_D} (\Omega) &= \{v\in H^1(\Omega);\; v|_{\Gamma_D} = g_D\},\\ 
\bH_{N,0}(\mathrm{div},\Omega) &= \{\bq\in \bH(\mathrm{div},\Omega);\; \langle\bq\cdot \bn ,v\rangle = 0\; \forall v \in H^1_{D,0}(\Omega)\},\\
\bH_{N,g_N}(\mathrm{div},\Omega)& = \{\bq\in \bH(\mathrm{div},\Omega);\; \langle\bq\cdot \bn ,v\rangle = g_N\; \forall v \in H^1_{D,0}(\Omega)\},
\end{align*}
and $\langle\cdot ,\cdot\rangle$ denotes the duality between $H^{-1/2}(\Gamma)$ and $H^{1/2}(\Gamma)$.
Under proper regularity assumptions, the above mixed formations are equivalent to each other, and each one can be adopted to construct Mixed DPGM. Notice that the boundary conditions are naturally built-in Mixed Formulation 4, and there is no need to calculate the derivatives for $u$ and $\bp$, which is an advantage for implementing the M-DPGM. Now, let us consider Mixed Formulation 4 to implement M-DPGM, and other mixed formulations can be applied similarly.

The Mixed Formulation 4 can be rewritten as: Find $(\bp,u) \in (L^2(\Omega))^d\times L^2(\Omega)$ such that
\begin{align}
	\mathcal{L}\big((\bp,u);(\bq,v)\big) = \mathcal{F}(\bq,v)\qquad \forall (\bq,v)\in \bH_{N,0}(\mathrm{div},\Omega)\times H^1_{D,0}(\Omega).\label{prb9}
\end{align}
where 
\begin{align*}
	\mathcal{L}\big((\bp,u);(\bq,v)\big) &= \int_{\Omega}(\bp\cdot \bq + u\,\nabla\cdot \bq + \bp\cdot\nabla v)\ {\rm d}\bx, \\
	\mathcal{F}(\bq,v) &= \int_{\Omega} f\,v\ {\rm d}\bx + \int_{\Gamma_N} g_N\,v\ \ds + \int_{\Gamma_D} g_D\,\bq\cdot \bn\ \ds.
\end{align*}

Similar to the DPGM introduced in Section \ref{prem}, we need to find two neural networks to approximate the variables $u$ and $\bp$, separately. Therefore,
the M-DPGM is: Find neural networks $u_\rho \in U_\rho\subset \mathcal{N}_{\rho}(D,M_D,B_{D})$ and $\bp_\rho \in \bQ_\rho\subset \mathcal{N}_{\rho}(\tilde{D},N_{\tilde{D}},B_{\tilde{D}},)$ such that
\begin{align}
	\mathcal{L}\big((\bp_\rho,u_\rho);(\bq_i,v_k)\big) = \mathcal{F}(\bq_i,v_k)\qquad \forall (\bq_i,v_k)\in \bQ_h\times V_h,\label{prb10}
\end{align}
where $\bQ_h$ and $V_h$ can be chosen as any proper finite-dimensional function spaces.
We construct $u_\rho \in U_\rho$ same as the one in \eqref{nnappu}, i.e., 
$$u_\rho(\bx) = \sum\limits_{j=1}^{ n_{D-1}} u_\rho^j \Phi^u_j(\bx).$$
Let us construct $\bp_\rho$ slightly different from $u_\rho$. Specifically, the number of neurons in the last layer for $\bp_\rho$ is $n_{\tilde{D}} = d$. Here, let us consider 2-dimensional case, i.e., $n_{\tilde{D}} = 2$, then the output of the neural network for $\bp_\rho$ is
\begin{equation*}
	\mathbf{\Phi}^p:= \mathbf{\Phi}^p_{\tilde{D}}(\mathbf{\bx}) =\mathbf{W}_{\tilde{D}}\mathbf{\Phi}^p_{{\tilde{D}}-1} = 
	\left[\begin{array}{c}\sum\limits_{j=1}^{ n_{{\tilde{D}}-1}}p_{1\rho}^{j} \Phi_{j}^{p}(\bx) \\
	\sum\limits_{j=1}^{ n_{{\tilde{D}}-1}}p_{2\rho}^{j} \Phi_{j}^{p}(\bx)\end{array}\right].
\end{equation*}
Here, we denote $w_{1j}^{D}$ and $w_{2j}^{D}$ by $p_{1\rho}^{j}$ and $p_{2\rho}^{j}$, respectively. Therefore, the problem \eqref{prb10} becomes: Find $u_\rho^j $ with $j = 1,\cdots,n_{D-1}$, and $p_{1\rho}^{j}$, $p_{2\rho}^{j}$ with $j = 1,\cdots,n_{{\tilde{D}}-1}$ such that
\begin{align}\label{prb11}
	\mathcal{L}\left(\left(\left[\sum\limits_{j=1}^{ n_{{\tilde{D}}-1}}p_{1\rho}^{j} \Phi_{j}^{p}(\bx),
	\sum\limits_{j=1}^{ n_{{\tilde{D}}-1}}p_{2\rho}^{j} \Phi_{j}^{p}(\bx)\right]^T,\sum\limits_{j=1}^{ n_{D-1}} u_\rho^j \Phi^u_j(\bx)\right);(\bq_i,v_k)\right) = \mathcal{F}(\bq_i,v_k)\qquad \forall (\bq_i,v_k)\in \bQ_h\times V_h.
\end{align}
Denote $\bq_i = [q_i^1,q_i^2]^T$, and we can take $(\bq_i,v_k) \in \bQ_h\times V_h$ in the forms of $([q^1_i,0],0),([0,q^2_i],0),([0,0],v_k)$, separately. Of course, other forms of the test functions can be used as well. Finally, we need to solve a least-square problem with the linear system generated by \eqref{prb11}.
\begin{remark}
\begin{enumerate}
\item Note that different neural network gives a different group of bases, if one wants to approximate several unknown functions in the same function space, one neural network can be used to approximate all the unknown functions at the same time, but the number of neurons in the last layer needs to be changed accordingly. 

\item For mixed finite element methods, the finite element pair $\bQ_h\times V_h$ needs to be carefully chosen so that the discrete problem is well-posed. Because the resulted linear system is solved by the least-square method, M-DPGM does not need to worry about the inf-sup condition, so it is very flexible on the choice of $\bQ_h\times V_h$.
\end{enumerate}
\end{remark}

\section{DPGM for time-dependent problems}
\label{sec:time}
In this section, we extend the ideas of DPGM to solve time-dependent PDEs. Instead of doing temporal and spatial discretization separately, we apply DPGM to solve time-dependent problems under the space-time approach, so the initial conditions will be treated as boundary conditions of the space-time domain. 

Let us consider the following time-dependent problem
\begin{align}
	&\mathcal{A} u = f \quad \;\;\;{\rm in}\; \Omega\times I, \label{pde3}\\
	&\mathcal{B} u = g \qquad {\rm on}\; \Gamma\times I, \label{pde4} \\
	&\mathcal{C} u = h \quad\;\; \;\;{\rm in}\; \Omega\times \{0\},  \label{pde5}
\end{align}
where $\Omega\subset \mathbb{R}^d$ is a bounded domain with $\Gamma = \Gamma_D \cup \Gamma_N $ and ${\Gamma}_D \cap {\Gamma}_N = \emptyset$, and $I=(0,T)$ is the time interval of interest. For example, when the differential operator $\mathcal{A}$ is given by
\begin{equation}
	\mathcal{A} :=\frac{\partial}{\partial t} -\nabla\cdot(\alpha(\bx)\nabla),
\end{equation}
\eqref{pde3} is a heat equation. Here, $\alpha_1\geq\alpha(\bx)\geq \alpha_0>0$ with some constants $\alpha_0$, $\alpha_1$. The boundary condition \eqref{pde4} can be Dirichlet, Neumann and Robin types, and the initial condition \eqref{pde5} is given by
\begin{align}
u(\bx,0) = h_0(\bx) \qquad {\rm in} \; \Omega.
\end{align}

Unlike the traditional approach, for example, finite difference method, discretizing the time interval as $0=t_0<t_1<\cdots<t_N=T$ and sequentially solving $u(\bx; t_n)$ for $n=1, \cdots, N$, we adopt space-time approach, i.e., temporal and spatial variables are treated jointly and equally.

\subsection{DPGM for a heat equation}

First, we consider to solve a heat equation as follows
\begin{align}
\frac{\partial u}{\partial t} -\nabla\cdot(\alpha(\bx)\nabla u ) &= f \quad\quad\;\;\;\;\;\; {\rm in}\; \Omega\times I, \\
u(\bx,t) &=  g_D(\bx,t)   \; \;\; {\rm on}\; \Gamma_D\times I,\\
\alpha(\bx)\nabla u(\bx,t)\cdot \bn &= g_N(\bx,t) \;\;\; {\rm on}\;\Gamma_N\times I,\\
u(\bx,0) &= h_0(\bx) \;\;\;\; \; \;\;{\rm in}\;  \Omega.
\end{align}
The weak formulation of the problem is to find $u\in L^2(I; H^1_{D,g_D}(\Omega))$ with $\partial_t u \in L^2(I;L^2(\Omega))$ such that
\begin{align}\label{heat_weak}
a_t(u,v) &= l_t(v) \qquad \forall v \in L^2(I;H^1_{D,0}(\Omega)),\\
u(\bx,0) &= h_0(\bx) \;\;\;\; \; \;{\rm in}\;  \Omega,
\end{align}
where
\begin{align*}
a_t(u,v)&=\int_{0}^{T}\int_\Omega \left(\frac{\partial u}{\partial t}v + \alpha(\bx)\nabla u\cdot \nabla v \right)\, {\rm d}\bx\, {\rm d} t,\\
l_t(v) &=\int_{0}^{T}\int_\Omega f\, v\,{\rm d}\bx\,{\rm d} t + \int_{0}^{T}\int_{\Gamma_N} g_N\, v\,\ds\,{\rm d} t.
\end{align*} 
Here, $L^2(I;H^1_{D,g_D}(\Omega))$ is a Banach space of all measurable functions $v: I\rightarrow H^1_{D,g_D}$ such that 
$$
\left[\int_0^T \|v(t)\|_{H^1(\Omega)}^2\; dt\right]^{1/2} < \infty.
$$

The DPGM for solving the problem \eqref{heat_weak} is to find a neural network $u_\rho(\bx,t) \in U_\rho\subset \mathcal{N}_{\rho}(D,M_D,B_{D}) \subset L^2(I;H^1(\Omega))$ such that 
\begin{align}
a_t(u_\rho,v) &= l_t(v) \qquad \quad\; \forall v \in V_h \subset L^2(I;H^1_{D,0}(\Omega)),\label{DGPM_heat_1}\\
u_\rho(\bx_k,t_k) &= g_D(\bx_k,t_k) \;\;\, {\rm for\ some\ points}\ (\bx_k,t_k)\in \Gamma_D  \times I,\; k=1,2,\cdots, N_b,\,\label{DGPM_heat_2}\\
u_{\rho}(\bx_m,0) &= h_0(\bx_m) \quad \;\;\,\,{\rm for\ some\ points}\ \bx_m\in  \Omega.\; m=1,2,\cdots,N_c,\, \label{DGPM_heat_3}
\end{align} 
Let us construct $u_\rho \in U_\rho$ same as the one in \eqref{nnappu} except for the input dimension $n_0 = d+1$ , i.e., 
$$u_\rho(\bx,t) = \sum\limits_{j=1}^{ n_{D-1}} u_\rho^j \Phi^u_j(\bx,t).$$
Similarly, the equation \eqref{DGPM_heat_1} becomes: Find $u_\rho^j $, $j = 1,\cdots,n_{D-1}$, such that
\begin{equation*}
\sum\limits_{j=1}^{ n_{D-1}}u_\rho^j  a_t(\Phi_j^u(\bx,t),v_i) = l_t(v_i)\qquad \forall v_i\in V_h,
\end{equation*}
which gives 
$$\mathbf{A} \boldsymbol{U}=\mathbf{L},$$ 
where $\mathbf{A}_{i,j} = a_t(\Phi_j^u(\bx,t),v_i)$, $\boldsymbol{U}= (u_\rho^1,\cdots,u_\rho^{n_{D-1}})^{T}$, and $\mathbf{L}=(l_1,\cdots,l_{N_h})^{T}$ with $l_i =l_t(v_i)$ ($i=1,\cdots,N_h$). Note that $\{\Phi_j^u(\bx,t)\}_{j=1}^{n_{D-1}}$ can be regarded as a group of  bases in space-time domain.

To deal with the boundary condition \eqref{DGPM_heat_2} and the initial condition \eqref{DGPM_heat_3}, we still need to take some random samples $\{(\bx_k,t_k)\}_{k=1}^{N_b}$ and $\{\bx_m\}_{m=1}^{N_c}$ according to the uniform distribution $\mathcal{U}(\Gamma_D \times I)$ and $\mathcal{U}(\Omega)$, separately. Then set $u_\rho(\bx_k,t_k)  =g_D(\bx_k,t_k)$ for all $ (\bx_k,t_k) \in \Gamma_D \times I$ and $u_\rho(\bx_m,0)  =h_0(\bx_m)$ for all $ \bx_m \in \Omega$, i.e.,
\begin{align*}
\sum\limits_{j=1}^{ n_{D-1}}u_\rho^j \Phi_j^u(\bx_k,t_k)&= g_D(\bx_k,t_k)  \qquad \text{for} \  k = 1,\cdots,N_b,\\
\sum\limits_{j=1}^{ n_{D-1}}u_\rho^j \Phi_j^u(\bx_m,0)&= h_0(\bx_m)  \qquad \;\;\;\;\text{for} \  m = 1,\cdots,N_c.
\end{align*}
Thus, the above equations imply $\mathbf{B} \boldsymbol{U}=\mathbf{G}$ and $\mathbf{C} \boldsymbol{U}=\mathbf{H}$, where $\mathbf{B}_{k,j} = \Phi_j^u(\bx_k,t_k)$, $\mathbf{C}_{m,j} = \Phi_j^u(\bx_m,0)$, $\mathbf{G}=(g_1,\cdots,g_{N_b})^{T}$ with  $g_k  =g_D(\bx_k,t_k)$ and $\mathbf{H}=(h_1,\cdots,h_{N_c})^{T}$ with  $h_m =h_0(\bx_m)$.

Finally we can obtain the solution $u_\rho$ by solving a least-square problem with the linear system 
\begin{equation}
\left[\begin{array}{c}\mathbf{A} \\\mathbf{B}\\\mathbf{C} \end{array}\right] \boldsymbol{U}
=\left[\begin{array}{c}\mathbf{L} \\\mathbf{G} \\\mathbf{H}\end{array}\right]. \label{equ}
\end{equation}

\subsection{DPGM for a wave equation}

In this subsection, we consider the following wave equation
\begin{align}
	\frac{\partial^2 u}{\partial t^2} -\nabla\cdot(\alpha(\bx)\nabla u ) &= f \quad\;\;\;\;\;\;\;\,\;\;\; {\rm in}\; \Omega\times I, \\
	 u(\bx,t) &= g_D(\bx,t)\quad\, {\rm on}\;\Gamma_D\times I, \\
	 \alpha(\bx)\nabla u(\bx,t)\cdot \bn &= g_N(\bx,t) \;\;\;\; {\rm on}\;\Gamma_N\times I,\\
	 	u(\bx,0) &= h_0(\bx)\qquad \,{\rm in}\;  \Omega,  \\
	 	\frac{\partial u(\bx,0)}{\partial t} &= w_0(\bx)\quad\;\;\;\;{\rm in}\;  \Omega. \label{ini_vel}
                 \end{align}
The weak form of the wave equation is: Find $u\in L^2(I; H^1_{D,g_D}(\Omega))$ with $\partial_t u \in L^2(I;L^2(\Omega))$ such that
\begin{align}
a_w(u,v) &= l_w(v) \qquad \forall v  \in L^2(I;H^1_{D,0}(\Omega)),\label{weak_wave}\\
u(\bx,0) &= h_0(\bx)\qquad {\rm in}\;  \Omega,  
\end{align}
where
\begin{align*}
a_w(u,v)&=\int_{0}^{T}\int_\Omega \left(-\frac{\partial u}{\partial t}\frac{\partial v}{\partial t} + \alpha(\bx)\nabla u\cdot \nabla v \right)
\, {\rm d}\bx\, {\rm d}t  + \int_\Omega \frac{\partial u(\bx,T)}{\partial t} v(\bx,T) \, {\rm d}\bx,\\
l_w(v) &=\int_{0}^{T}\int_\Omega f\, v\,{\rm d}\bx\,{\rm d} t+ \int_{0}^{T}\int_{\Gamma_N} g_N\, v\,\ds\,{\rm d}t  + \int_\Omega w_0v(\bx,0) \, {\rm d}\bx.
\end{align*} 
Note that the initial condition \eqref{ini_vel} becomes a natural boundary condition of the domain $\Omega\times I$, and it is built into the weak formulation \eqref{weak_wave}.

The DPGM for solving the wave equation is to find a neural network $u_\rho \in U_\rho\subset \mathcal{N}_{\rho}(D,M_D,B_{D})\subset L^2(I;H^1(\Omega))$ such that 
\begin{align*}
a_w(u_\rho,v) &= l_w(v) \qquad \;\;\;\forall v  \in V_h \subset L^2\left(I;H^1_{D,0}(\Omega)\right),\\
u_\rho(\bx_k,t_k) &= g_D(\bx_k,t_k) \;\;\, {\rm for\ some\ points}\ (\bx_k,t_k)\in \Gamma_D  \times I,\; k=1,2,\cdots, N_b,\,\\
u_{\rho}(\bx_m,0) &= h_0(\bx_m) \quad \;\;\,\,{\rm for\ some\ points}\ \bx_m\in  \Omega.\; m=1,2,\cdots,N_c,\, 
\end{align*}
By constructing $u_{\rho}(\bx,t)$ and making restrictions on the boundary and initial conditions in a similar way for the heat equation, we obtain a similar linear system as \eqref{equ}, then solve the least-square problem to get solution $u_\rho$. 

\begin{remark}
Usually, one prefers dividing a big problem into some small problems and solving them one by one so that memory cost is low and total computation time is saved. For example, finite difference discretization for the temporal variable is the most popular choice for the time-dependent problem because it follows this divide and conquers strategy. However, one drawback is the accumulation of errors during this step-by-step process. If we can solve the problem at just one time, this issue can be avoided, but large computation ability is needed to support this approach, for example, the space-time FEM is used to solve the time-dependent problem with the help of domain decomposition. In light of the strong approximation ability of neural networks, DPGM can solve time-dependent problems by the space-time approach efficiently and accurately.
\end{remark}

\section{Numerical examples}

In this section, we present three examples, a Poisson equation, a heat equation, and a wave equation for demonstrating the performance of the DPGM. In these examples, we consider $d=2$ and choose basis functions of finite element method as the test functions, specifically, the domain $\Omega$ and $\Omega\times I$ are decomposed into the square and cubic elements, respectively, and corresponding bilinear functions and trilinear functions are adopted as test functions.
Other types of test functions can be considered as well in the future.

For the calculation of $\nabla u_\rho$, one can adopt the automatic differentiation, which allows the calculation of derivatives for a broad range of functions. However, in order to obtain the fast computation of derivatives, we employ the difference method, for example, the central difference for a first-order derivative is
\begin{equation*}
\frac{\partial{u_\rho(x,y)}}{\partial x} = \frac{u_\rho(x+h_x,y)-u_\rho(x-h_x,y)}{2h_x}.
\end{equation*}
In numerical experiments, we find that it saves lots of computation time with less loss of accuracy by choosing $h_x =h_y =10^{-6}$.

Initializing the network with the right weights is crucial for us, once the wight parameters are initialized, they will be fixed in DPGM except for the last layer. We must make sure that these weights are in a reasonable range before the solving process. In practical experiments, one may choose proper initialization from $torch.nn.init$ in Pytorch. In the following examples, we use the uniform distribution for the fully connected network, and Xavier uniform distribution is chosen for ResNet based on its property of remaining variance the same during each passing layer (\cite{glorot2010xavier}). The activation function $\rho = tanh$ is used in the fully connected network and ResNet as it leads to a smooth $u_{\rho}$ for approximating the solution of PDEs.

\begin{example}\label{exam1}In this example, we solve a 2-dimensional Poisson equation with a smooth solution $u = cos(\pi x)sin(\pi y)$ under the mixed boundary conditions
	\begin{numcases}{}
		-\Delta u(x,y) = 2\pi^2cos(\pi x)sin(\pi y) \quad {\rm in} \; \Omega, \notag \\
		u(x,y) =cos(\pi x)sin(\pi y) \qquad\qquad {\rm on} \;\Gamma_D, \notag\\ 
		\nabla u(x,y)\cdot \bn = 0 \qquad\qquad\qquad\quad\;\;\;{\rm on} \;\Gamma_N,\notag
	\end{numcases}
	where $\Omega = [0,1]^2$, $\Gamma_D = (0,1) \times \{0,1\}$, $\Gamma_N = \{0,1\}\times (0,1) $.
\end{example}

For DPGM, we use a two-layer fully connected neural network with uniform distribution $\mathcal{U}(-1, 1)$ as initialization, $D = 2$, $n_0 = 2$, $n_2 = 1$, and we choose $n_{1} = 50,100,200$, separately. We know that the degrees of freedom (dof) for the DPGM is the number of neurons in the first layer, that is, $n_{1}$. Test functions are chosen as bases of bilinear finite element on a family of square meshes with mesh size $h=2^{-n}\; (n=2,3,4,5)$, thus the result will be related to $h$. We calculate the numerical integration by Gauss-Legendre quadrature with 25 points inside of each square, and randomly sample 100 points on each edge of $\Gamma_D$. For solving the resulted linear system, we adopt the least-square solver $scipy.linalg.lstsq$ in Python.

\begin{figure}[!ht]
	\begin{center}
 	\subfigure[Numerical solution $u_\rho$]{
		 \label{bar_domain}
		 \centering
		 \includegraphics[width=2.0in]{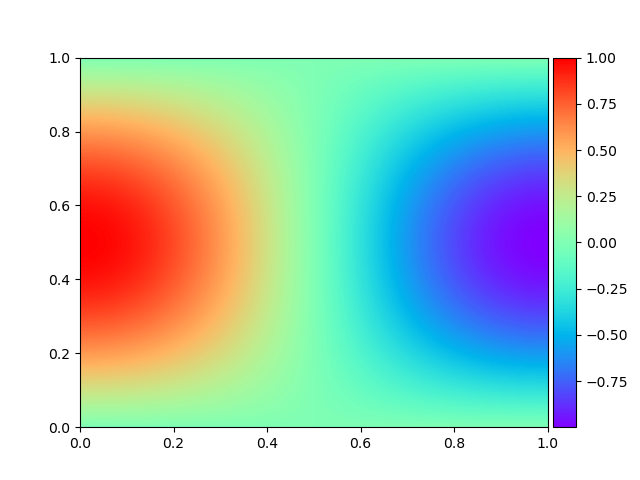}
  	}
 	\subfigure[Exact solution $u$]{
	 	\label{vuggy_domain}
	 	\centering
	 	\includegraphics[width=2.0in]{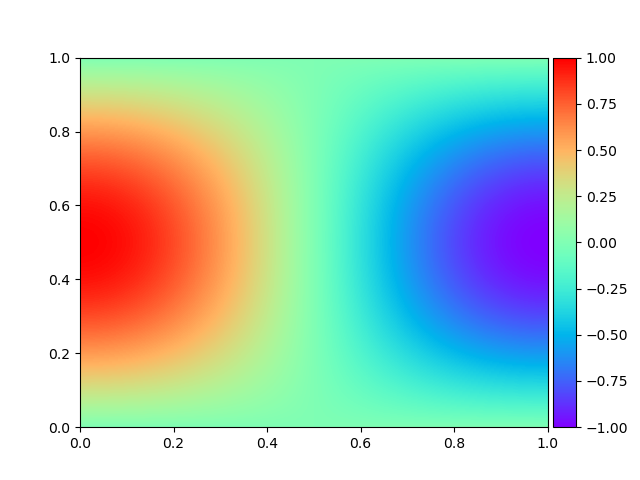}
 	}
	\subfigure[$|u-u_\rho|$]{
	 	\label{vuggy_domain}
	 	\centering
	 	\includegraphics[width=2.0in]{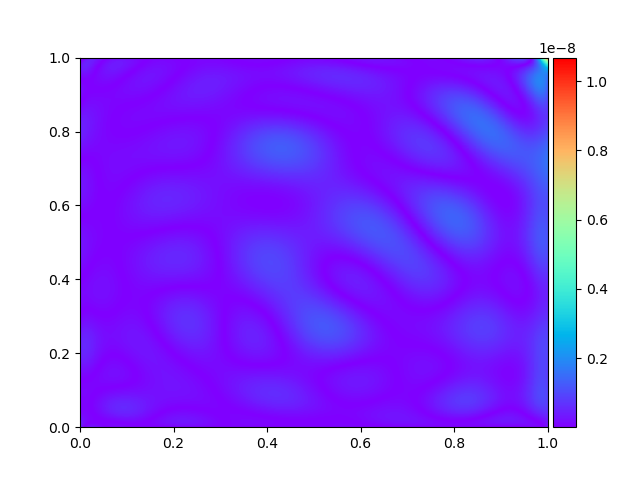}
 	}
	\end{center}
	\vspace*{-15pt}
	\caption{Numerical solution $u_{\rho}$ by DPGM with $h=2^{-5}$ and dof $=200$ in Example \ref{exam1}.}
	\label{figure1}
\end{figure}

To demonstrate the accuracy of the DPGM, we give the numerical solution $u_\rho$, the exact solution $u$, and their difference $|u-u_\rho|$ in \rm{Figure} \ref{figure1}. The Figure \ref{figure1} (a) shows the numerical solution $u_{\rho}$ solved by the DPGM with $h = 2^{-5}$ and dof $= 200$, and its relative $L^2$ error and $H^1$ error are $9.651\times10^{-10}$ and $5.167\times10^{-9}$, respectively. The exact solution $u$ is given in Figure \ref{figure1} (b), and the absolute value of their difference $|u-u_\rho|$ is displayed in Figure \ref{figure1} (c), which shows that the maximum error is around $10^{-8}$. We see that the DPGM gives highly accurate numerical solution. 

In addition, in Table \ref{table1a}, we report the relative $L^2$ and $H^1$ errors denoted by $e_{L^2}$ and $e_{H^1}$, respectively, with $h=2^{-n}\; (n=2,3,4,5)$ and dof $=50,100,200$. With fixed $h$, i.e., given the same data of test functions to the neural network, DPGM offers more accurate numerical solutions with a larger number of dof = $n_{1}$. If we fix the number of dof, and change the mesh size $h$ for the test functions, we observe that the errors decrease as $h$ becomes smaller, that is, more test functions are fed into the system, the neural networks supply more accurate numerical solutions.

Of course, we want to know the performance of this new DPGM compared with the finite element method, so we compute some numerical solutions by FEM through Fenics (\cite{langtangena2017fenics}). On a family of uniform triangulation with mesh size $h=2^{-n}\; (n=2,3,4,5)$, we use the standard triangle Lagrange elements $P_k$ ($k=1,2,3$) for solving the problem. Here, $k$ denotes the degree of the polynomials, and the degrees of freedom for $P_k$ FEM in this example is dof$ = (k/h)^2-1$). The numerical errors are listed in Table \ref{table1b}. In view of \rm{Table} \ref{table1a} and \rm{Table} \ref{table1b}, comparing the dof and accuracy, we observe that the DPGM outperforms over FEM, that is, DPGM can obtain more accurate numerical solution with much less degrees of freedom.


\begin{table}[!htbp]	
	\centering  
	\begin{tabular}{|c|c|c|c|c|c|c|}  
		\hline  
		\diagbox [width=5em] {$h$}{dof}&  
		\multicolumn{2}{c|}{50}&\multicolumn{2}{c|}{100}&
		\multicolumn{2}{c|}{200}\cr\cline{1-7}  
		&$e_{L^2}$ & $e_{H^1}$ &$e_{L^2}$ & $e_{H^1}$&$e_{L^2}$ & $e_{H^1}$ \cr\hline  
		$2^{-2}$&6.895e-3 &2.080e-2&4.745e-3 &1.507e-2&2.931e-3 &9.062e-3 \cr\hline  
		$2^{-3}$&4.604e-5 &3.161e-4&1.116e-5 &7.537e-5&4.081e-6 &2.732e-5\cr\hline  
		$2^{-4}$&2.892e-5 &1.139e-4&2.684e-8 &1.219e-7&2.890e-9 &1.383e-8 \cr\hline  
		$2^{-5}$&2.730e-5&9.634e-5&1.219e-8 &5.417e-8&9.651e-10 &5.167e-9\cr\hline  
	\end{tabular}  
	\caption{Relative $L^2(\Omega)$ and $H^1(\Omega)$ errors of the DPGM with different $h$ and dof in Example \ref{exam1}.}
	\label{table1a}
\end{table}

\begin{table}[!htbp] 	
	\centering  
	\begin{tabular}{|c|c|c|c|c|c|c|c|c|c|}  
		\hline  
		\diagbox [width=6em] {$h$}{Scheme}&  
		\multicolumn{3}{c|}{$P_1$ Linear FEM}&\multicolumn{3}{c|}{$P_2$ quadratic FEM}&
		\multicolumn{3}{c|}{$P_3$ cubic FEM}\cr\cline{1-10}  
		&dof&$e_{L^2}$ & $e_{H^1}$ &dof&$e_{L^2}$ & $e_{H^1}$&dof&$e_{L^2}$ & $e_{H^1}$ \cr\hline  
		$2^{-2}$&15& 1.829e-1 &3.990e-1&63 &8.878e-3 &5.615e-2&143&6.626e-4 &5.724e-3 \cr\hline  
		$2^{-3}$&63&4.399e-2 &1.934e-1&255&1.102e-3 & 1.456e-2&575&3.958e-5 & 7.202e-4 \cr\hline  
		$2^{-4}$&255&1.091e-2 &9.601e-2 &1023&1.375e-4 &3.683e-3&2303&2.417e-6 & 9.004e-5 \cr\hline  
		$2^{-5}$&1023&2.721e-3 &4.792e-2&4095& 1.718e-5& 9.246e-4&9215&1.495e-7&1.125e-5 \cr\hline  
	\end{tabular}  
	\caption{Relative $L^2(\Omega)$ and $H^1(\Omega)$ errors of the $P_k$ FEM with different $h$ in Example \ref{exam1}.}
	\label{table1b}
\end{table}

Furthermore, we study the relationship between data and networks in the DPGM. Specifically, data is determined by the test functions, changing the mesh size or choosing different numbers of the test functions will generate different data, while the network can be formed by different depth and width, which means the number of neurons in each year.
Now, let us fix the mesh size as $h=2^{-5}$, and report the relative errors $L^2(\Omega)$ and $H^1(\Omega)$ for different number of unfixed neurons (dof) and the numbers of the test functions (nv) in \rm{Table} \ref{table1c}. We see that the more data information is supplied, the more accurate solution the DPGM will offer. Then we adopt the ResNet with Xavier uniform distribution $\widetilde{\mathcal{U}}(-1,1)$ as initialization and change the number of layers $D$ from 2 to 5, report the relative errors $L^2(\Omega)$ and $H^1(\Omega)$ for different neurons (dof) in each layer in \rm{Table} \ref{table1d}, and we find that the deeper neural network would give us a better numerical solution.

\begin{table}[!htbp]	
	\centering  
	\begin{tabular}{|c|c|c|c|c|c|c|}  
		\hline  
		\diagbox [width=5em] {$nv$}{dof}&  
		\multicolumn{2}{c|}{50}&\multicolumn{2}{c|}{100}&
		\multicolumn{2}{c|}{200}\cr\cline{1-7}  
		&$e_{L^2}$ & $e_{H^1}$ &$e_{L^2}$ & $e_{H^1}$&$e_{L^2}$ & $e_{H^1}$ \cr\hline  
		200&1.274e-4 &1.861e-4&3.043e-7 &8.408e-7&1.364e-7 &1.570e-7 \cr\hline  
		400&1.244e-4 &2.093e-4&1.449e-7 &3.028e-7&1.435e-8 &2.177e-8\cr\hline  
		600&1.005e-4 &2.881e-4&8.909e-8 &1.573e-7&3.745e-9 &1.387e-8 \cr\hline  
		800&7.368e-5 &2.020e-4&4.091e-8 &1.751e-7&3.968e-9 &1.065e-8\cr\hline  
		1023&3.914e-5 &1.437e-4&2.597e-8 &1.039e-7&9.651e-10 &5.167e-9\cr\hline  
	\end{tabular}  
	\caption{Relative $L^2(\Omega)$ and $H^1(\Omega)$ errors of the DPGM with respect to dof and the number of test functions (nv) for a fixed mesh size $h=2^{-5}$ in Example \ref{exam1}.}
	\label{table1c}
\end{table}

\begin{table}[!htbp]	
	\centering  
	\begin{tabular}{|c|c|c|c|c|c|c|}  
		\hline  
		\diagbox [width=5em] {$D$}{dof}&  
		\multicolumn{2}{c|}{50}&\multicolumn{2}{c|}{100}&
		\multicolumn{2}{c|}{200}\cr\cline{1-7}  
		&$e_{L^2}$ & $e_{H^1}$ &$e_{L^2}$ & $e_{H^1}$&$e_{L^2}$ & $e_{H^1}$ \cr\hline  
		2&1.178e-5 &4.100e-5 &5.827e-6 &1.405e-5 &9.318e-5 &2.054e-4
\cr\hline  
		3&1.093e-5 &3.165e-5 &1.661e-6 &4.765e-6&1.342e-6&3.937e-6
\cr\hline  
		4&1.302e-5 &4.786e-5 &1.857e-7 &7.485e-7 &
1.240e-7 &4.744e-7\cr\hline  
		5&6.548e-5 &2.674e-4 &4.485e-8 &2.093e-7
&4.629e-8 &2.108e-7\cr\hline  
		
	\end{tabular}  
	\caption{Relative $L^2(\Omega)$ and $H^1(\Omega)$ errors of the DPGM with respect to dof and the depth $D$ for a fixed mesh size $h=2^{-5}$ in Example \ref{exam1}.}
	\label{table1d}
\end{table}

Finally, under the same condition, we consider the mixed DPGM with different formulations, which show nice performance in \rm{Table} \ref{table1k}. We can see that all four schemes get better results with higher degrees of freedom. The relative $L^2$ error of M-DPGM-4 can even reduce to $1.395\times10^{-10}$, which shows its high accuracy, and note that M-DPGM-4 do not need the calculation of the derivatives for $u$ and $\bp$, and boundary conditions are naturally built in the formulation, so no boundary restriction is needed.

\begin{table}[!htbp] 	
	\centering  
	\begin{tabular}{|c|c|c|c|c|c|c|c|c|c|c|}  
		\hline  
		\diagbox [width=6em] {$h$}{Scheme}& &
		\multicolumn{2}{c|}{M-DPGM-1}&\multicolumn{2}{c|}{M-DPGM-2}&
		\multicolumn{2}{c|}{M-DPGM-3}& \multicolumn{2}{c|}{M-DPGM-4}\cr\cline{1-10}  
		&dof&$e_{L^2}$ & $e_{H^1}$ &$e_{L^2}$& $e_{H^1}$ & $e_{L^2}$&$e_{H^1}$& $e_{L^2}$ & $e_{H^1}$ \cr\hline  
		$2^{-2}$&75& 4.408e-2 &1.315e-1&7.098e-1 &1.677e-1&3.785e-3 &1.754e-2&7.333e-2 &7.569e-2 \cr\hline  
 		$2^{-3}$&150& 4.857e-5&2.063e-4&8.316e-5 &2.287e-4&3.103e-5 &9.329e-5&5.396e-5 &1.268e-4 \cr\hline  
		$2^{-4}$&300& 9.944e-8 &2.984e-7&5.140e-8&4.180e-7&4.989e-8 &3.502e-7 &1.995e-8 &1.004e-7 \cr\hline  
		$2^{-5}$&600& 3.503e-9 &7.141e-8&2.783e-9 &3.795e-8&4.086e-9 &4.536e-8&1.395e-10 &1.678e-9 \cr\hline  
	
	\end{tabular}  
	\caption{Relative $L^2(\Omega)$ and $H^1(\Omega)$ errors of different formulation of Mixed DPGM with different $h$ and dof in Example \ref{exam1}.}
	\label{table1k}
\end{table}

\begin{example}\label{exam2} Given $\Omega = (0,1)^2$ with $\Gamma_D = \partial \Omega:= \Gamma_1\cup \Gamma_2\cup \Gamma_3\cup \Gamma_4$, we consider the the following heat equation
	\begin{numcases}{}
		u_t(x,y,t)-\Delta u(x,y,t) = f(x,y,t) \quad {\rm in}\; \Omega\times I, \notag \\
		u(x,y,t) = g(x,y,t) \qquad \qquad \qquad \quad{\rm on}\;\Gamma_D\times I,\notag\\
		u(x,y,0) = h_0(x,y)\qquad\qquad \qquad\quad\,{\rm in}\; \Omega, \notag
	\end{numcases}
	with $I = (0,1)$.
	The exact solution $u(x,y,t) = 2e^{-t}sin(\frac{\pi}{2}x)sin(\frac{\pi}{2}y)$ and the right-hand term $f(x,y,t)$ is given accordingly.
\end{example}

As we mentioned in Section \ref{sec:time}, we do not discretize the temporal variable by finite difference method, instead, DPGM solves the above equation directly in the space-time domain, that is, we treat this 2-dimensional heat equation as a 3-dimensional problem with 2-dimensional spatial variables and 1-dimensional temporal variable. In the test, we use a two-layer fully connected neural network with uniform distribution $\mathcal{U}(-1, 1)$ as initialization, $D = 2$, $n_0 = 3$, $n_2 = 1$, and we choose $n_{1} = 200,400,800$, separately. The test functions are chosen by trilinear functions on cubic meshes and numerical integration is calculated by Gauss-Legendre quadrature with 1000 points inside of each cubic. Note that this space-time domain has 6 faces, to enforce the Dirichlet and initial condition, we randomly sample 100 points on each face, i.e., $\Gamma_1\times I$, $\Gamma_2\times I$, $\Gamma_3\times I$, $\Gamma_4\times I$, and $\Omega\times\{0\}$. The $L^2$ errors and $H^1$ errors at the ending time $T=1$ are reported in Table \ref{table2a}. We observe that the DPGM works very well for solving heat equation, and still reach high accuracy. 

To compare the performance of DPGM with traditional numerical methods, we consider a numerical scheme that is constructed with the $P_2$ finite element discretization for the spatial variable and the back Euler finite difference approximation for the temporal variable.
We choose the time-steps $\Delta t = 10^{-3}, 2\times 10^{-4}, 5\times 10^{-5}$ and the mesh-size $h = 2^{-5}, 2^{-6}, 2^{-7}$, then report the $L^2$ errors and $H^1$ errors at the ending time $T=1$. Compared to FEM, DPGM can obtain more accurate numerical solutions, furthermore, the time cost of DPGM is much less than the traditional approaches.

\begin{table}[!htbp]	
	\centering  
	\begin{tabular}{|c|c|c|c|c|c|c|}  
		\hline  
		\diagbox [width=6em] {$h$}{dof}&  
		\multicolumn{2}{c|}{200}&\multicolumn{2}{c|}{400}&
		\multicolumn{2}{c|}{800}\cr\cline{1-7}  
		&$L^2$ error & $H^1$ error &$L^2$ error & $H^1$ error &$L^2$ error & $H^1$ error \cr\hline  
		$2^{-2}$&9.577e-5&1.368e-3&5.184e-5 &7.290e-4&3.448e-5&4.952e-4 \cr\hline  
		$2^{-3}$&2.605e-5&2.871e-4&1.274e-7 &1.967e-6&6.490e-9 &1.673e-7\cr\hline  
		$2^{-4}$&2.845e-5 &2.600e-4&1.743e-7&2.004e-6&8.347e-10 &1.544e-8 \cr\hline  
	\end{tabular}  
	\caption{$L^2(\Omega)$ and $H^1(\Omega)$ errors of the DPGM with different $h$ and dof in Example \ref{exam2}.}
	\label{table2a}
\end{table}

\begin{table}[!htbp]	
	\centering  
	\begin{tabular}{|c|c|c|c|c|c|c|}  
		\hline  
		\diagbox [width=6em] {$h$}{$\Delta t$}&  
		\multicolumn{2}{c|}{$10^{-3}$}&\multicolumn{2}{c|}{$2\times 10^{-4}$}&
		\multicolumn{2}{c|}{$5\times 10^{-5}$}\cr\cline{1-7}  
		&$L^2$ error & $H^1$ error &$L^2$ error & $H^1$ error & $L^2$ error & $H^1$ error \cr\hline  
		$2^{-5}$&7.268e-6 &3.481e-5&1.450e-6 &7.209e-6&3.597e-7&2.657e-6 \cr\hline  
		$2^{-6}$&7.272e-6 &3.480e-5&1.454e-6 & 6.960e-6&3.632e-7 &1.756e-6 \cr\hline 
		$2^{-7}$&7.272e-6 &3.480e-5&1.454e-6 & 6.958e-6& 3.635e-7& 1.739e-6 \cr\hline 
	\end{tabular}  
	\caption{$L^2(\Omega)$ and $H^1(\Omega)$ errors of the $P_2$ FEM for spatial discretization with different $h$ and back Euler scheme for temporal discretization with different time-step $\Delta t$ in Example \ref{exam2}}
	\label{table2b}
\end{table}

\begin{example}\label{exam3} Consider the following wave equation 
\begin{numcases}{}
	\frac{\partial^2 u}{\partial t^2}(x,y,t)-\Delta u(x,y,t) = f(x,y,t) \;\;\; {\rm in}\; \Omega\times (0,1), \notag \\
	u(x,y,t) = g(x,y,t) \qquad \qquad \qquad \quad\;\;\,{\rm on}\;\Gamma_D\times (0,1),\notag\\
	u(x,y,0) = h_0(x,y)\qquad\qquad \qquad\qquad{\rm in}\; \Omega, \notag\\
	\frac{\partial u}{\partial t}(x,y,0) = w_0(x,y)\qquad\qquad \qquad\;\;\;\;{\rm in}\; \Omega, \notag
\end{numcases}
with the exact solution $u(x,y,t) = sin(\frac{\pi}{2}x)sin(\frac{\pi}{2}y)sin(\frac{\pi}{2}t)$	over $\Omega = (0,1)^2 \subset \mathbb{R}^2$ and $\Gamma_D = \partial \Omega$. The right-hand term $f(x,y,t)$ is given accordingly. 
\end{example}

In this example, we apply the DPGM based on the space-time approach to solve the wave equation. we use a two-layer fully connected neural network with uniform distribution $\mathcal{U}(-1, 1)$ as initialization, $D = 2$, $n_0 = 3$, $n_2 = 1$, and we choose $n_{1} = 200,400,800$, separately. Same as the case of the heat equation in Example \ref{exam2}, the test functions are chosen as trilinear functions on cubic meshes. Same as Example 5.2, we randomly sample 100 points on each faces enforcing the Dirichlet boundary and the initial conditions. The $L^2$ and $H^1$ errors at the ending time $T=1$ are shown in Table \ref{table3a}, and we can see that the DPGM still obtain very accurate numerical solutions.

\begin{table}[!htbp]	
	\centering  
	\begin{tabular}{|c|c|c|c|c|c|c|}  
		\hline  
		\diagbox [width=6em] {$h$}{$dof$}&  
		\multicolumn{2}{c|}{200}&\multicolumn{2}{c|}{400}&
		\multicolumn{2}{c|}{800}\cr\cline{1-7}  
		&$L^2$ error & $H^1$ error &$L^2$ error & $H^1$ error &$L^2$ error & $H^1$ error \cr\hline  
		$2^{-2}$&5.818e-5 &8.190e-4&5.461e-5 &7.700e-4&4.672e-5 &6.525e-4 \cr\hline  
		$2^{-3}$&1.015e-5 &1.085e-4&6.175e-8 &1.760e-6&4.796e-9 &1.446e-7\cr\hline  
	 $2^{-4}$&1.839e-5& 2.412e-4&5.632e-8 &1.018e-6&7.290e-10 &2.495e-8\cr\hline  
	\end{tabular}  
	\caption{$L^2(\Omega)$ and $H^1(\Omega)$ errors of the DPGM with different $h$ and dof in Example \ref{exam3}.}
	\label{table3a}
\end{table}

\section{Summary}

This new framework, the deep Petrov-Galerkin method, uses neural networks to approximate solutions of partial differential equations, then solve the resulting linear system by the least-square method. The DPGM is based on variational formulation, the trial function space is approximated by neural networks and test function space can be flexibly given by other numerical methods or neural networks. The resulted linear system is not symmetric and square, so the discretized problem has to be solved by the least-square method. Compared with other numerical methods, for example, finite element method, finite difference method, DPGM has the following advantages: (i) it supplies much more accurate numerical solution with respect to degrees of freedom due to the powerful approximation property of neural networks; (ii) this method is mesh-free and the choice of test functions is flexible; (iii) boundary conditions can be treated easily in this framework; (iv) mixed DPGM can be easily constructed to approximate several unknown functions simultaneously, and it has good flexibility to handle different boundary conditions; (v) DPGM can solve the time-dependent problems by space-time approach naturally and efficiently. 

We believe that this new numerical framework has a strong potential for solving various partial differential equations, however, this newborn baby is immature, and needs to be taken care of carefully in many aspects. What is the performance if other neural networks and test functions are used? The resulted linear system may have a large condition number, how do we design neural networks and choose proper test functions to avoid this situation? We know that there are many finite-dimensional function spaces used in the spectral method, like Fourier basis functions for periodic problems, Chebyshev or Legendre polynomials for problems defined in bounded domains, Lagrange polynomials for problems in semi-bounded domains, Hermite polynomials for problems in an unbounded domain, could we use these basis functions in the DPGM? The numerical analysis of this method is quite open and needs to be explored further. How to design adaptive DPGM to improve its efficiency? What is the performance of this method for solving other more complex problems? Can we implement this method in parallel for solving large-scale problems? Many works related to this method are waiting for us to explore.

\end{document}